\documentclass{amsart}
\usepackage{amsfonts}

\newtheorem{theorem}{Theorem}
\theoremstyle{plain}

\newtheorem{corollary}{Corollary}

\newtheorem{definition}{Definition}
\newtheorem{example}{Example}

\newtheorem{lemma}{Lemma}

\newtheorem{remark}{Remark}

\numberwithin{equation}{section}
\input{tcilatex}

\begin{document}
\title{Characterizations of $\Gamma $-AG$^{\ast \ast }$-groupoids by their $%
\Gamma $-ideals}
\author{}
\maketitle

\begin{center}
\textbf{Madad Khan, Naveed Ahmad and Inayatur Rehman}

\bigskip

COMSATS Institute of Information Technology

Abbottabad Pakistan

madadmath@yahoo.com, naveed.maths@yahoo.com

\bigskip 

Department of Mathematics, Quaid-i-Azam University

Islamabad, Pakistan

s\_inayat@yahoo.com

\bigskip
\end{center}

\textbf{Abstract.} In this paper we have discusses $\Gamma $-left, $\Gamma $%
-right, $\Gamma $-bi-, $\Gamma $-quasi-, $\Gamma $-interior and $\Gamma $%
-ideals in $\Gamma $-AG$^{\ast \ast }$-groupoids and regular $\Gamma $-AG$%
^{\ast \ast }$-groupoids. Moreover we have proved that the set of $\Gamma $%
-ideals in a regular $\Gamma $-AG$^{\ast \ast }$-groupoid form a semilattice
structure. Also we have characterized a regular $\Gamma $-AG$^{\ast \ast }$%
-groupoid in terms of left ideals.

\section{Introduction}

Kazim and Naseeruddin \cite{Kaz} have introduced the concept of an
LA-semigroup. This structure is the generalization of a commutative
semigroup. It is closely related with a commutative semigroup and
commutative groups because if an LA-semigroup contains right identity then
it becomes a commutative semigroup and if a new binary operation is defined
on a commutative group which gives an LA-semigroup \cite{Mushtaq and Yusuf}.
The connection of the class of LA-semigroups with the class of vector spaces
over finite fields and fields has been given as: Let $W$ be a sub-space of a
vector space $V$ over a field $F$ of cardinal $2r$ such that $r>1$. Many
authors have generalized some useful results of semigroup theory.

In 1981, the notion of $\Gamma $-semigroups was introduced by M. K. Sen \cite%
{gemma1} and \cite{gamma2}.

T. Shah and I. Rehman \cite{irehman} defined $\Gamma $-AG-groupoids
analogous to $\Gamma $-semigroups and then they introduce the notion of $%
\Gamma $-ideals and $\Gamma $-bi-ideals in $\Gamma $-AG-groupoids. It is
easy to see that $\Gamma $-ideals and $\Gamma $-bi-ideals in $\Gamma $%
-AG-groupoids are infect a generalization of ideals and bi-ideals in
AG-groupoids (for a suitable choice of $\Gamma $).

In this paper we define $\Gamma $-quasi-ideals and $\Gamma $-interior ideals
in $\Gamma $-AG$^{\ast \ast }$-groupoids and generalize some results. Also
we have proved that $\Gamma $-AG-groupoids with left identity and
AG-groupoids with left identity coincide.

Let $G$ and $\Gamma $ be two non-empty sets. $G$ is said to be a $\Gamma $%
-AG-groupoid if there exist a mapping $G\times \Gamma \times G\rightarrow G$%
, written $\left( a\text{, }\gamma \text{, }b\right) $ as $a\gamma b$, such
that $G$ satisfies the identity $\left( a\gamma b\right) \delta c=\left(
c\gamma b\right) \delta a$, for all $a$, $b$, $c\in G$ and $\gamma $, $%
\delta \in \Gamma $ \cite{irehman}.

\begin{definition}
An element $e\in S$ is called a left identity of $\Gamma $-AG-groupoid if $%
e\gamma a=a$ for all $a\in S$ and $\gamma \in \Gamma $.
\end{definition}

\begin{lemma}
If a $\Gamma $-AG-groupoid contains left identity, then it becomes an
AG-groupoid with left identity.
\end{lemma}

\begin{proof}
Let $G$ be a $\Gamma $-AG-groupoid and $e$ be the left identity of $G$ and
let $a$, $b\in G$ and $\alpha $, $\beta \in \Gamma $ therefore we have%
\begin{equation*}
a\alpha b=a\alpha (e\beta b)=e\alpha (a\beta b)=a\beta b\text{.}
\end{equation*}

Hence $\Gamma $-AG-groupoid with left identity becomes and an AG-groupoid
with left identity.
\end{proof}

\begin{remark}
From Lemma 1, it is easy to see that all the results given in \cite{irehman}
and \cite{Rehman2} for a $\Gamma $-AG-groupoid with left identity is
identical to the results given in \cite{Madad1} and \cite{Madad2}.
\end{remark}

\begin{definition}
A $\Gamma $-AG-groupoid is called a $\Gamma $-AG$^{\ast \ast }$-groupoid if
it satisfies the following law%
\begin{equation*}
a\alpha (b\beta c)=b\alpha (a\beta c)\text{, for all }a,b,c\in S\text{ and }%
\alpha ,\beta \in \Gamma \text{.}
\end{equation*}
\end{definition}

The following results and definition from definition \ref{1} to lemma \ref%
{idempotent} have been taken from \cite{irehman}.

\begin{definition}
\label{1}Let $G$ be a $\Gamma $-AG-groupoid, a non-empty subset $S$ of $G$
is called sub $\Gamma $-AG-groupoid if $a\gamma b\in S$ for all $a$, $b\in S$
and $\gamma \in \Gamma $ or $S$ is called sub $\Gamma $-AG-groupoid if$\
S\Gamma S\subseteq S.$
\end{definition}

\begin{definition}
A subset $I$ of a $\Gamma $-AG-groupoid $G$ is called left(right) $\Gamma $%
-ideal of $G$ if $G\Gamma I\subseteq I\left( I\Gamma G\subseteq I\right) $
and $I$ is called $\Gamma $-ideal of $G$ if it is both left and right $%
\Gamma $-ideal.
\end{definition}

\begin{definition}
An element $a$ of a $\Gamma $-AG-groupoid $G$ is called regular if there
exist $x\in G$ and $\beta $, $\gamma \in \Gamma $ such that $a=\left( a\beta
x\right) \gamma a$. $G$ is called regular $\Gamma $-AG-groupoid if all
elements of $G$ are regular.
\end{definition}

\begin{definition}
A sub $\Gamma $-AG-groupoid $B$ of a $\Gamma $-AG-groupoid $G$ is called $%
\Gamma $-bi-ideal of $G$ if $\left( B\Gamma G\right) \Gamma B\subseteq B$.
\end{definition}

\begin{definition}
Let $G$ and $\Gamma $ be any non-empty sets. If there exists a mapping $%
G\times \Gamma \times G\rightarrow G,$ written $\left( x\text{, }\gamma 
\text{, }y\right) $ as $x\gamma y$, $G$ is called a $\Gamma $-medial if it
satisfies $\left( x\alpha y\right) \beta \left( l\gamma m\right) =\left(
x\alpha l\right) \beta \left( y\gamma m\right) $, and called $\Gamma $%
-paramedial if it satisfies $\left( x\alpha y\right) \beta \left( l\gamma
m\right) =\left( m\alpha l\right) \beta \left( y\gamma x\right) $ for all $x$%
, $y$, $l$, $m\in G$ and $\alpha $, $\beta $, $\gamma \in \Gamma $.
\end{definition}

\begin{lemma}
\label{ab=ba}If $A$ and $B$ are any $\Gamma $-ideals of a regular $\Gamma $%
-AG-groupoid $G$ then $A\Gamma B=B\Gamma A$.
\end{lemma}

\begin{definition}
A $\Gamma $-ideal $P$ of a $\Gamma $-AG-groupoid $G$ is called $\Gamma $%
-prime$\left( \Gamma \text{-semiprime}\right) $ if for any $\Gamma $-ideals $%
A$ and $B$, $A\Gamma B\subseteq P\left( A\Gamma A\subseteq P\right) $
implies either $A\subseteq P$ or $B\subseteq P\left( A\subseteq P\right) $.
\end{definition}

\begin{lemma}
\label{idempotent}Any $\Gamma $-ideal $A$ of a regular $\Gamma $-AG-groupoid
is a $\Gamma $-idempotent that is $A\Gamma A=A$.
\end{lemma}

It is important to note that every $\Gamma $-AG-groupoid $G$ is $\Gamma $%
-medial and every $\Gamma $-AG$^{\ast \ast }$-groupoid $G$ is $\Gamma $%
-paramedial because for any $x$, $y$, $l$, $m\in G$ and $\alpha $, $\beta $, 
$\gamma \in \Gamma $, we have 
\begin{equation*}
\left( x\alpha y\right) \beta \left( l\gamma m\right) =\left( \left( l\gamma
m\right) \alpha y\right) \beta x=\left( \left( y\gamma m\right) \alpha
l\right) \beta x=\left( x\alpha l\right) \beta \left( y\gamma m\right) \text{%
.}
\end{equation*}

We call it as $\Gamma $-medial law.

\begin{theorem}
If $L$ and $R$ are left and right $\Gamma $-ideals of a $\Gamma $-AG$^{\ast
\ast }$-groupoid $G$ then $L\cup L\Gamma G$ and $R\cup G\Gamma R$ are $%
\Gamma $-ideals of $G$.
\end{theorem}

\begin{proof}
Let $L$ be a left $\Gamma $-ideal of $G$ then we have%
\begin{eqnarray*}
\left( L\cup L\Gamma G\right) \Gamma G &=&\left( L\Gamma G\right) \cup
\left( L\Gamma G\right) \Gamma G=\left( L\Gamma G)\cup (G\Gamma G\right)
\Gamma L \\
&\subseteq &L\Gamma G\cup \left( G\Gamma L\right) \subseteq L\Gamma G\cup
L=L\cup L\Gamma G\text{ and} \\
G\Gamma \left( L\cup L\Gamma G\right) &=&G\Gamma L\cup G\Gamma \left(
L\Gamma G\right) \subseteq L\cup L\Gamma \left( G\Gamma G\right) =L\cup
L\Gamma G.
\end{eqnarray*}

Again let $R$ be a right $\Gamma $-ideal of $G$ then we have%
\begin{eqnarray*}
\left( R\cup G\Gamma R\right) \Gamma G &=&R\Gamma G\cup \left( G\Gamma
R\right) \Gamma G\subseteq R\cup \left( G\Gamma R\right) \Gamma \left(
G\Gamma G\right) \\
&=&R\cup \left( G\Gamma G\right) \Gamma \left( R\Gamma G\right) \subseteq
R\cup G\Gamma R\text{, and } \\
G\Gamma \left( R\cup G\Gamma R\right) &=&G\Gamma R\cup G\Gamma \left(
G\Gamma R\right) =G\Gamma R\cup \left( G\Gamma G\right) \Gamma \left(
G\Gamma R\right) \\
&=&G\Gamma R\cup \left( R\Gamma G\right) \Gamma \left( G\Gamma G\right)
\subseteq G\Gamma R\cup R\Gamma G \\
&\subseteq &G\Gamma R\cup R=R\cup G\Gamma R\text{.}
\end{eqnarray*}
\end{proof}

\begin{lemma}
\label{right}Right identity in a $\Gamma $-AG-groupoid $G$\ becomes identity
of $G$ and hence $G$ becomes commutative $\Gamma $-semigroup.
\end{lemma}

\begin{proof}
Let $e$ be the right identity of $G$ , $g\in G$, $\alpha $ and $\beta \in
\Gamma $, then 
\begin{equation*}
e\alpha g=\left( e\beta e\right) \alpha g=\left( g\beta e\right) \alpha
e=g\alpha e=g.
\end{equation*}

Again for $a$, $b$, $c\in G$ and $\alpha $, $\beta \in \Gamma $ we have%
\begin{equation*}
a\gamma b=\left( e\alpha a\right) \gamma b=\left( e\alpha a\right) \gamma
(e\alpha b)=(b\alpha e)\gamma \left( a\alpha e\right) =b\gamma a\text{.}
\end{equation*}

Now 
\begin{eqnarray*}
\left( a\alpha b\right) \beta c &=&\left( a\alpha b\right) \beta \left(
e\alpha c\right) =\left( a\alpha e\right) \beta \left( b\alpha c\right)
=e\alpha \left( \left( a\alpha e\right) \beta \left( b\alpha c\right) \right)
\\
&=&\left( a\alpha e\right) \alpha \left( e\beta \left( b\alpha c\right)
\right) =a\alpha \left( e\beta \left( b\alpha c\right) \right) =a\alpha
\left( b\beta \left( e\alpha c\right) \right) \\
&=&a\alpha \left( b\beta c\right) \text{.}
\end{eqnarray*}
\end{proof}

\begin{definition}
A sub $\Gamma $-AG-groupoid $Q$ of a $\Gamma $-AG-groupoid $G$ is called a
quasi-ideal of $G$ if $G\Gamma Q\cap Q\Gamma G\subseteq Q$.
\end{definition}

\begin{definition}
A sub $\Gamma $-AG-groupoid $I$ of a $\Gamma $-AG-groupoid $G$ is called a $%
\Gamma $-interior ideal of $G$ if $\left( G\Gamma I\right) \Gamma G\subseteq
I$.
\end{definition}

\begin{lemma}
Every one sided (left or right) $\Gamma $-ideal of a $\Gamma $-AG-groupoid $%
G $ is a $\Gamma $-quasi ideal of $G$.
\end{lemma}

\begin{proof}
Let $L$ be a left $\Gamma $-ideal of $G$ then we have%
\begin{equation*}
L\Gamma G\cap G\Gamma L\subseteq G\Gamma L\subseteq L\text{.}
\end{equation*}

Which implies $L$ is a $\Gamma $-quasi ideal of $G$. Similarly if $R$ is a
right $\Gamma $-ideal of $G\,$then it is a $\Gamma $-quasi ideal of $G$.
\end{proof}

\begin{lemma}
\label{RLB}Every right $\Gamma $-ideal and left $\Gamma $-ideal\ of a $%
\Gamma $-AG-groupoid $G$ is a $\Gamma $-bi-ideal of $G$.
\end{lemma}

\begin{proof}
Let $R$ be a right $\Gamma $-ideal of $G$ then we have%
\begin{equation*}
\left( R\Gamma G\right) \Gamma R\subseteq R\Gamma R\subseteq R\Gamma
G\subseteq R\text{.}
\end{equation*}

Again let $L$ be a left $\Gamma $-ideal of $G$ then we have%
\begin{equation*}
\left( L\Gamma G\right) \Gamma L\subseteq \left( G\Gamma G\right) \Gamma
L\subseteq G\Gamma L\subseteq L\text{.}
\end{equation*}
\end{proof}

\begin{corollary}
Every $\Gamma $-ideal of a $\Gamma $-AG-groupoid $G$ is a $\Gamma $-bi-ideal
of $G$.
\end{corollary}

\begin{proof}
It follows from lemma \ref{RLB}.
\end{proof}

\begin{lemma}
If $B_{1}$ and $B_{2}$ are $\Gamma $-bi-ideals of a $\Gamma $-AG$^{\ast \ast
}$-groupoid $G$ then $B_{1}\Gamma B_{2}$ is also a $\Gamma $-bi-ideals of $G$%
.
\end{lemma}

\begin{proof}
Let $B_{1}$ and $B_{2}$ be $\Gamma $-bi-ideals of $G$ then we have 
\begin{eqnarray*}
\left( \left( B_{1}\Gamma B_{2}\right) \Gamma G\right) \Gamma \left(
B_{1}\Gamma B_{2}\right) &=&\left( \left( B_{1}\Gamma B_{2}\right) \Gamma
\left( G\Gamma G\right) \right) \Gamma \left( B_{1}\Gamma B_{2}\right) \\
&=&\left( \left( B_{1}\Gamma G\right) \Gamma \left( B_{2}\Gamma G\right)
\right) \Gamma \left( B_{1}\Gamma B_{2}\right) \\
&=&\left( \left( B_{1}\Gamma G\right) \Gamma B_{1}\right) \Gamma \left(
\left( B_{2}\Gamma G\right) \Gamma B_{2}\right) \\
&\subseteq &B_{1}\Gamma B_{2}.
\end{eqnarray*}
\end{proof}

\begin{lemma}
Every $\Gamma $-idempotent quasi-ideal of a $\Gamma $-AG-groupoid $G$ is a $%
\Gamma $-bi-ideal of $G$.
\end{lemma}

\begin{proof}
Let $Q$ be an $\Gamma $-idempotent quasi-ideal of $G$. Now%
\begin{eqnarray*}
\left( Q\Gamma G\right) \Gamma Q &\subseteq &\left( G\Gamma G\right) \Gamma
Q\subseteq G\Gamma Q\text{, and} \\
\left( Q\Gamma G\right) \Gamma Q &=&\left( Q\Gamma G\right) \Gamma \left(
Q\Gamma Q\right) =\left( Q\Gamma Q\right) \Gamma \left( G\Gamma Q\right)
=Q\Gamma \left( G\Gamma Q\right) \\
&\subseteq &Q\Gamma \left( G\Gamma G\right) \subseteq Q\Gamma G\text{, which
implies that } \\
\left( Q\Gamma G\right) \Gamma Q &\subseteq &G\Gamma Q\cap Q\Gamma
G\subseteq Q\text{.}
\end{eqnarray*}
\end{proof}

\begin{lemma}
\label{ideal interior}Every $\Gamma $-ideal of a $\Gamma $-AG-groupoid $G$
is a $\Gamma $-interior ideal of $G$.
\end{lemma}

\begin{proof}
Let $I$ be a $\Gamma $-ideal of $G$ then we have%
\begin{equation*}
\left( G\Gamma I\right) \Gamma G\subseteq I\Gamma G=I.
\end{equation*}
\end{proof}

\begin{lemma}
A subset $I$ of a $\Gamma $-AG$^{\ast \ast }$-groupoid $G$ is a $\Gamma $%
-interior ideal if and only if it is right $\Gamma $-ideal.
\end{lemma}

\begin{proof}
Let $I$ be a right $\Gamma $-ideal $G$ then it becomes a left $\Gamma $%
-ideal so is $\Gamma $-ideal and by lemma \ref{ideal interior} it is $\Gamma 
$-interior ideal.

Conversely assume that $I$ is a $\Gamma $-interior ideal of $G$. Using $%
\Gamma $-paramedial law, we have 
\begin{eqnarray*}
I\Gamma G &=&I\Gamma \left( G\Gamma G\right) =G\Gamma \left( I\Gamma
G\right) =\left( G\Gamma G\right) \Gamma \left( I\Gamma G\right) \\
&=&\left( G\Gamma I\right) \Gamma \left( G\Gamma G\right) \subseteq \left(
G\Gamma I\right) \Gamma G\subseteq G.
\end{eqnarray*}

Which shows that $I$ is a right $\Gamma $-ideal of $G$.
\end{proof}

\begin{example}
Let $G=\left\{ 1,2,3,4,5\right\} $ with binary operation "$\cdot $" given in
the following Cayley's table, an AG-groupoid with left identity $4$.%
\begin{equation*}
\begin{tabular}{l|lllll}
$\cdot $ & $1$ & $2$ & $3$ & $4$ & $5$ \\ \hline
$1$ & $4$ & $5$ & $1$ & $2$ & $3$ \\ 
$2$ & $3$ & $4$ & $5$ & $1$ & $2$ \\ 
$3$ & $2$ & $3$ & $4$ & $5$ & $1$ \\ 
$4$ & $1$ & $2$ & $3$ & $4$ & $5$ \\ 
$5$ & $5$ & $1$ & $2$ & $3$ & $4$%
\end{tabular}%
\end{equation*}
\end{example}

It is easy to observe that $G$ is a simple AG-groupoid that is there is no
left or right ideal of $G$. Now let $\Gamma =\{\alpha $, $\beta $, $\gamma
\} $ defined as 
\begin{equation*}
\begin{tabular}{l|lllll}
$\alpha $ & $1$ & $2$ & $3$ & $4$ & $5$ \\ \hline
$1$ & $1$ & $1$ & $1$ & $1$ & $1$ \\ 
$2$ & $1$ & $1$ & $1$ & $1$ & $1$ \\ 
$3$ & $1$ & $1$ & $1$ & $1$ & $1$ \\ 
$4$ & $1$ & $1$ & $1$ & $1$ & $1$ \\ 
$5$ & $1$ & $1$ & $1$ & $1$ & $1$%
\end{tabular}%
\text{ \ \ \ \ \ \ \ \ }%
\begin{tabular}{l|lllll}
$\beta $ & $1$ & $2$ & $3$ & $4$ & $5$ \\ \hline
$1$ & $2$ & $2$ & $2$ & $2$ & $2$ \\ 
$2$ & $2$ & $2$ & $2$ & $2$ & $2$ \\ 
$3$ & $2$ & $2$ & $2$ & $2$ & $2$ \\ 
$4$ & $2$ & $2$ & $2$ & $2$ & $2$ \\ 
$5$ & $2$ & $2$ & $2$ & $2$ & $2$%
\end{tabular}%
\text{ \ \ \ \ \ \ \ \ }%
\begin{tabular}{l|lllll}
$\gamma $ & $1$ & $2$ & $3$ & $4$ & $5$ \\ \hline
$1$ & $1$ & $1$ & $1$ & $1$ & $1$ \\ 
$2$ & $1$ & $1$ & $1$ & $1$ & $1$ \\ 
$3$ & $1$ & $1$ & $1$ & $1$ & $1$ \\ 
$4$ & $1$ & $1$ & $1$ & $1$ & $1$ \\ 
$5$ & $1$ & $1$ & $1$ & $3$ & $3$%
\end{tabular}%
\end{equation*}

It is easy to prove that $G$ is a $\Gamma $-AG-groupoid because $\left( a\pi
b\right) \psi c=\left( c\pi b\right) \psi a\ $for all $a$, $b$, $c\in G$ and 
$\pi $, $\psi \in \Gamma $ also $G$ is non-associative because $\left(
1\alpha 2\right) \beta 3\not=1\alpha \left( 2\beta 3\right) $. This $\Gamma $%
-AG-groupoid does not contain left identity because $4\alpha 5\neq 5$, $%
4\beta 5\neq 5$ and $4\gamma 5\neq 5$. It is easy to see that every
AG-groupoid with left identity not necessarily implies $\Gamma $%
-AG-groupoid\ with left identity. Clearly $A=\{1,2,3\}$ is a $\Gamma $-ideal
of $G$. $B=\{1,2,4\}$ is a right $\Gamma $-ideal but is not a left $\Gamma $%
-ideal. $A$ and $B$ both are $\Gamma $-bi-ideals of $G$. $C=\{1$, $2$, $3$, $%
4\}$ is a $\Gamma $-interior ideal of $G$.

\begin{lemma}
For a regular $\Gamma $-AG-groupoid $G$ $A\Gamma G=A$ and $G\Gamma B=B$ for
every right $\Gamma $-ideal $A$ and for every left $\Gamma $-ideal $B$.
\end{lemma}

\begin{proof}
Let $A$ be a right $\Gamma $- ideal of $G$ then $A\Gamma G\subseteq A$. Let $%
a\in A$, since $G$ is regular so there exist $x\in G$ and $\alpha $, $\gamma
\in \Gamma $ such that%
\begin{equation*}
a=\left( a\alpha x\right) \gamma a\in \left( A\Gamma G\right) \Gamma
A\subseteq \left( A\Gamma G\right) \Gamma G\subseteq A\Gamma G\text{.}
\end{equation*}

Now again let $B$ be a left $\Gamma $-ideal of $G$ then $G\Gamma B\subseteq
B $. Let $b\in B$, also $G$ is regular so there exist $t\in G$ and $\pi $, $%
\sigma \in \Gamma $ such that 
\begin{equation*}
b=\left( b\pi t\right) \sigma b\in \left( B\Gamma G\right) \Gamma B\subseteq
\left( G\Gamma G\right) \Gamma B\subseteq G\Gamma B\text{.}
\end{equation*}
\end{proof}

\begin{lemma}
If $G$ is a $\Gamma $-AG$^{\ast \ast }$-groupoid then $g\Gamma G$ and $%
G\Gamma g$ are $\Gamma $-bi-ideals for all $g\in G$.
\end{lemma}

\begin{proof}
Using the definition of $\Gamma $-AG$^{\ast \ast }$-groupoid we have 
\begin{eqnarray*}
\left( \left( g\Gamma G\right) \Gamma G\right) \Gamma \left( g\Gamma
G\right) &=&\left( \left( G\Gamma G\right) \Gamma g\right) \Gamma \left(
g\Gamma G\right) \subseteq \left( G\Gamma g\right) \Gamma \left( g\Gamma
G\right) \\
&=&g\Gamma \left( \left( G\Gamma g\right) \Gamma G\right) \subseteq g\Gamma
\left( \left( G\Gamma G\right) \Gamma G\right) \subseteq g\Gamma \left(
G\Gamma G\right) \\
&\subseteq &g\Gamma G\text{.}
\end{eqnarray*}

Again using $\Gamma $-paramedial law we have%
\begin{eqnarray*}
\left( \left( G\Gamma g\right) \Gamma G\right) \Gamma \left( G\Gamma
g\right) &=&\left( \left( \left( G\Gamma g\right) \Gamma g\right) \Gamma
G\right) \Gamma G=\left( \left( \left( g\Gamma g\right) \Gamma G\right)
\Gamma G\right) \Gamma G \\
&=&\left( \left( G\Gamma G\right) \Gamma G\right) \Gamma \left( g\Gamma
g\right) \subseteq \left( G\Gamma G\right) \Gamma \left( g\Gamma g\right) \\
&=&\left( g\Gamma g\right) \Gamma \left( G\Gamma G\right) \subseteq \left(
g\Gamma g\right) \Gamma G=\left( G\Gamma g\right) \Gamma g \\
&\subseteq &\left( G\Gamma G\right) \Gamma g\subseteq G\Gamma g.
\end{eqnarray*}
\end{proof}

\begin{corollary}
If $G$ is a regular $\Gamma $-AG$^{\ast \ast }$-groupoid then $a\Gamma G$ is
a $\Gamma $-bi-ideal in $G$, for all $a\in G$.
\end{corollary}

\begin{proof}
Let $G$ be a regular $\Gamma $-AG-groupoid then for every $a\in G$ there
exist $x\in G$ and $\alpha $, $\beta \in \Gamma $ such that $a=\left( \left(
a\alpha x\right) \beta a\right) $ therefore we have%
\begin{eqnarray*}
\left( \left( a\Gamma G\right) \Gamma G\right) \Gamma \left( a\Gamma
G\right) &=&\left( \left( \left( \left( a\alpha x\right) \beta a\right)
\Gamma G\right) \Gamma G\right) \Gamma \left( a\Gamma G\right) \\
&=&\left( \left( G\Gamma G\right) \Gamma \left( \left( a\alpha x\right)
\beta a\right) \right) \Gamma \left( a\Gamma G\right) \\
&\subseteq &\left( G\Gamma \left( \left( a\alpha x\right) \beta a\right)
\right) \Gamma \left( a\Gamma G\right) =\left( \left( a\alpha x\right)
\Gamma \left( G\beta a\right) \right) \Gamma \left( a\Gamma G\right) \\
&\subseteq &\left( \left( a\alpha x\right) \Gamma \left( G\beta G\right)
\right) \Gamma \left( G\Gamma G\right) \subseteq \left( \left( a\alpha
x\right) \Gamma G\right) \Gamma G \\
&=&\left( G\Gamma G\right) \Gamma \left( a\alpha x\right) \subseteq G\Gamma
\left( a\alpha x\right) =a\Gamma \left( G\alpha x\right) \subseteq a\Gamma
\left( G\Gamma G\right) \\
&\subseteq &a\Gamma G\text{.}
\end{eqnarray*}
\end{proof}

\begin{lemma}
For a $\Gamma $-bi-ideal $B$ in a regular $\Gamma $-AG-groupoid $G$, $\left(
B\Gamma G\right) \Gamma B=B$.
\end{lemma}

\begin{proof}
Let $B$ be a $\Gamma $-bi-ideal in $G$ then $\left( B\Gamma G\right) \Gamma
B\subseteq B$. Let $x\in B$, since $G$ is a regular $\Gamma $-AG-groupoid
therefore there exist $a\in G$ and $\alpha $, $\beta \in \Gamma $ such that 
\begin{equation*}
x=\left( x\alpha a\right) \beta x\in \left( B\Gamma G\right) \Gamma B\text{.}
\end{equation*}

Which implies that $B\subseteq \left( B\Gamma G\right) \Gamma B$.
\end{proof}

\begin{lemma}
If $G$ is a regular $\Gamma $-AG-groupoid then, $G\Gamma G=G$.
\end{lemma}

\begin{proof}
Since $G\Gamma G\subseteq G.$ Let $x\in G$, since $G$ is a regular $\Gamma $%
-AG-groupoid therefore there exist $a\in G$ and $\alpha $, $\beta \in \Gamma 
$ such that 
\begin{equation*}
x=\left( x\alpha a\right) \beta x\in \left( G\Gamma G\right) \Gamma
G\subseteq G\Gamma G\text{.}
\end{equation*}

Which implies that $G\subseteq G\Gamma G$.
\end{proof}

\begin{lemma}
A subset $I$ of a regular $\Gamma $-AG$^{\ast \ast }$-groupoid $G$ is a left 
$\Gamma $-ideal if and only if it is a right $\Gamma $-ideal of $G.$
\end{lemma}

\begin{proof}
Let $I$ be a left $\Gamma $-ideal of $G$ then $G\Gamma I\subseteq I$. Let $%
i\gamma g\in I\Gamma G$ for $g\in G$, $i\in I$ and $\gamma \in \Gamma $,
also $G$ is a regular $\Gamma $-AG-groupoid therefore there exist $x$, $y\in
G$ and $\alpha $, $\beta $, $\gamma $, $\delta $, $\pi \in \Gamma $ such that%
\begin{eqnarray*}
i\gamma g &=&\left( \left( i\alpha x\right) \beta i\right) \gamma \left(
\left( g\delta y\right) \pi g\right) =\left( \left( i\alpha x\right) \beta
\left( g\delta y\right) \right) \gamma \left( i\pi g\right) \\
&=&\left( \left( \left( \left( i\alpha x\right) \beta i\right) \alpha
x\right) \beta \left( g\delta y\right) \right) \gamma \left( i\pi g\right)
=\left( \left( y\alpha g\right) \beta \left( \left( i\beta \left( i\alpha
x\right) \right) \delta x\right) \right) \gamma \left( i\pi g\right) \\
&=&\left( i\beta \left( \left( \left( y\alpha g\right) \beta \left( i\alpha
x\right) \right) \delta x\right) \right) \gamma \left( i\pi g\right) =\left(
\left( i\pi g\right) \beta \left( \left( \left( y\alpha g\right) \beta
\left( i\alpha x\right) \right) \delta x\right) \right) \gamma i \\
&\in &\left( G\Gamma I\right) \subseteq I\text{.}
\end{eqnarray*}

Conversely let $I$ be a right $\Gamma $-ideal then there exist $x\in G$ and $%
\alpha $, $\beta \in \Gamma $ such that%
\begin{equation*}
g\gamma i=\left( \left( g\alpha x\right) \beta g\right) \gamma i=\left(
i\beta g\right) \gamma \left( g\alpha x\right) \in \left( I\Gamma G\right)
\Gamma G\subseteq I\Gamma G\subseteq I\text{.}
\end{equation*}
\end{proof}

\begin{theorem}
for a $\Gamma $-AG$^{\ast \ast }$-groupoid $G$, following statements are
equivalent.

$\left( i\right) $ $G$ is regular $\Gamma $-AG-groupoid.

$\left( ii\right) $ Every left $\Gamma $-ideal of $G$ is $\Gamma $-idempotent%
$.$
\end{theorem}

\begin{proof}
$\left( i\right) \Rightarrow \left( ii\right) $

Let $G$ be a regular $\Gamma $-AG-groupoid then by lemma \ref{idempotent}
every $\Gamma $-ideal of $G$ is $\Gamma $-idempotent.

$\left( ii\right) \Rightarrow \left( i\right) $

Let every left $\Gamma $-ideal of a $\Gamma $-AG$^{\ast \ast }$-groupoid $G$
is $\Gamma $-idempotent, since $G\Gamma a$ is a left $\Gamma $-ideal of $G$
for all $a\in G$ \cite{irehman}, so is $\Gamma $-idempotent and by $\Gamma $%
-paramedial law, lemma \ref{a(bc)} and $\Gamma $-medial law, we have, $a\in
G\Gamma a$ implies 
\begin{eqnarray*}
a &\in &\left( G\Gamma a\right) \Gamma \left( G\Gamma a\right) =\left(
\left( G\Gamma a\right) \Gamma a\right) \Gamma G=\left( \left( a\Gamma
a\right) \Gamma G\right) \Gamma G \\
&=&\left( \left( a\Gamma a\right) \Gamma \left( G\Gamma G\right) \right)
\Gamma G=\left( \left( G\Gamma G\right) \Gamma \left( a\Gamma a\right)
\right) \Gamma G \\
&=&\left( a\Gamma \left( \left( G\Gamma G\right) \Gamma a\right) \right)
\Gamma G=\left( G\Gamma \left( \left( G\Gamma G\right) \Gamma a\right)
\right) \Gamma a \\
&=&\left( G\Gamma \left( G\Gamma a\right) \right) \Gamma a=\left( G\Gamma
\left( \left( G\Gamma a\right) \Gamma \left( G\Gamma a\right) \right)
\right) \Gamma a \\
&=&\left( G\Gamma \left( \left( a\Gamma G\right) \Gamma \left( a\Gamma
G\right) \right) \right) \Gamma a=\left( \left( G\Gamma G\right) \Gamma
\left( \left( a\Gamma G\right) \Gamma \left( a\Gamma G\right) \right)
\right) \Gamma a \\
&=&\left( \left( G\Gamma \left( a\Gamma G\right) \right) \Gamma \left(
G\Gamma \left( a\Gamma G\right) \right) \right) \Gamma a=\left( \left(
\left( a\Gamma G\right) \Gamma G\right) \Gamma \left( \left( a\Gamma
G\right) \Gamma G\right) \right) \Gamma a \\
&=&\left( \left( \left( \left( a\Gamma G\right) \Gamma G\right) \Gamma
G\right) \Gamma \left( a\Gamma G\right) \right) \Gamma a=\left( a\Gamma
\left( \left( \left( \left( a\Gamma G\right) \Gamma G\right) \Gamma G\right)
\Gamma G\right) \right) \Gamma a \\
&\subseteq &\left( a\Gamma G\right) \Gamma a.
\end{eqnarray*}

Which shows that $G$ is a regular $\Gamma $-AG$^{\ast \ast }$-groupoid.
\end{proof}

\begin{lemma}
Any $\Gamma $-ideal of a regular $\Gamma $-AG-groupoid $G$ is $\Gamma $%
-semiprime.
\end{lemma}

\begin{proof}
It is an easy consequence of lemma \ref{idempotent}.
\end{proof}

\begin{theorem}
Set of all $\Gamma $-ideals in a regular $\Gamma $-AG-groupoid $G$ with
forms a semilattice $\left( G,\circ \right) $ where $A\circ B=A\Gamma B$,
for all $\Gamma $-ideals $A$ and $B$ of $G.$
\end{theorem}

\begin{proof}
Let $A$ and $B$ be any $\Gamma $-ideals in $G$, then by $\Gamma $-medial law
we have%
\begin{eqnarray*}
\left( A\Gamma B\right) \Gamma G &=&\left( A\Gamma B\right) \Gamma \left(
G\Gamma G\right) =\left( A\Gamma G\right) \Gamma \left( B\Gamma G\right)
\subseteq A\Gamma B.\text{ And } \\
G\Gamma \left( A\Gamma B\right) &=&\left( G\Gamma G\right) \Gamma \left(
A\Gamma B\right) =\left( G\Gamma A\right) \Gamma \left( G\Gamma B\right)
\subseteq A\Gamma B\text{.}
\end{eqnarray*}

Also by lemma \ref{ab=ba}, we have $A\Gamma B=B\Gamma A$ which implies that 
\begin{equation*}
\left( A\Gamma B\right) \Gamma C=C\Gamma \left( A\Gamma B\right) =A\Gamma
\left( C\Gamma B\right) =A\Gamma \left( B\Gamma C\right) \text{.}
\end{equation*}

And by lemma \ref{idempotent}, $A\Gamma A=A$.\newpage
\end{proof}

\end{document}